\documentstyle[12pt]{amsart}
\newcommand{\g}{\goth}

\newcommand{\gtsl}{\mbox{\g sl}}

\newcommand{\hgtsl}{\mbox{$\hat{\gtsl}$}}

\newcommand{\nc}{\mbox{${\Bbb C}$}}
\newcommand{\nz}{\mbox{${\Bbb Z}$}}

\newcommand{\cF}{\mbox{${\cal F}$}}

\newcommand{\cP}{\mbox{${\cal P}$}}

\newcommand{\vep}{\varepsilon}
\newcommand{\vphi}{\varphi}
\theoremstyle{plain}
 \newtheorem{thm}{Theorem}[section]
 \newtheorem{prop}[thm]{Proposition}
 \newtheorem{lemma}[thm]{Lemma}
 \newtheorem{cor}[thm]{Corollary}
\theoremstyle{definition}
 \newtheorem{defn}{Definition}[section]

\theoremstyle{remark}

\begin{document}
\title{Difference equations of quantum current operators and 
quantum parafermion construction
\\} 
\author{Jintai Ding}
\author{Boris Feigin}
\address{Jintai Ding, RIMS, Kyoto University}
\address{Boris Feigin, RIMS, Kyoto University; Landau Institute, Moscow}
\maketitle
\begin{abstract}
For the current realization of the affine quantum groups, a 
simple comultiplication for the quantum current operators was given 
by Drinfeld. With this 
comultiplication, we prove that, for the integrable modules of 
$U_q(\hat {\frak sl}(2))$ of level $k+1$, $x^\pm(z)x^\pm(zq^{\pm 2})
\cdot\cdot\cdot
x^\pm(zq^{\pm 2k})$ are vertex operators satisfying  certain q-difference 
equations, and we derive the quantum parafermions
of $U_q(\hat {\frak sl}(2))$.

\end{abstract}
\pagestyle{plain}
\section{Introduction.} 

Lie algebra $\hat {\frak sl}(2)$ has three current operators $e(z)$, 
$h(z)$ and $f(z)$. 
For any   integrable highest weight module of $\hat {\frak sl}(2)$ of
level $k$, the current operators  $e(z)$ and $f(z)$ satisfy the following 
differential equation:
$$(e^k(z))'-:z^{-1}h(z)e^k(z):=0=(f^k(z))'+:z^{-1}h(z)f(z): ,$$
which implies that $e^k(z)$ and $f^k(z)$ 
 are vertex operators \cite{LP}.

For the case of quantum affine algebras, Drinfeld 
presented a  formulation of affine quantum groups with generators 
in the form of current operators\cite{Dr2}, which, for the case of
  $U_q(\hat {\frak sl}(2))$,   give us the 
 quantized current operators corresponding to $e(z)$, $h(z)$ and $f(z)$ of 
 $\hat {\frak sl}(2)$. We would like to find out 
if it is possible to derive a similar 
equation, which will degenerate into the equation above. 
 To solve this problem, we need to use 
the Drinfeld comultiplication for the current formulation of 
$U_q(\hat {\frak sl}(2))$\cite{DF}\cite{DI},
 which resolve the difficulty caused by the non-commutativity of those 
quantum current operators. 
 This comultiplication is very simple
as opposed to  the comultiplication formula induced from the conventional
 comultiplication, which  
 can not be written in a closed form with those current operators.
With this comultiplication, we are able to 
study the zeros and poles of quantum current operators for 
integrable modules to derive a quantum integrable condition for 
  $U_q(\hat {\frak sl}(2))$\cite{DM}. 

In this paper, we will use the same method as in \cite{DM}. We will start 
with the case of the module $\c F$ for the  fundamental 
representations at the level 1 for the case of $U_q(\hat{\frak sl}(2))$, 
which is the representation constructed by Frenkel and Jing
by  using  vertex operators. $U_q(\hat{\frak sl}(2))$ as in the Drinfeld 
realization (see Definition 2.1) has four current generators 
$x^+(z)$, $\vphi(z)$, $\psi(z)$ and $x^-(z)$, where $x^+(z)$ and $x^-(w)$
 are quantized current operators of $U_q(\hat {\frak sl}(2))$
 corresponding to $e(z)$ and $f(z)$ of $\hat {\frak sl}(2)$ respectively
and  $\vphi(z)$, $\psi(z)$ are quantized current operator 
corresponding to the negative half and the positive half of 
$h(z)$ of $\hat {\frak sl}(2)$  respectively. 
Using Drinfeld comultiplication, we show that, on  
any level $m+1$ integrable module of 
$U_q(\hat {\frak sl}(2))$, 
$$(x^+(zq^{2}))(x^+(zq^{4}))\cdot\cdot\cdot
(x^+(zq^{(2m+2)}))
=$$
$$ \vphi^{-1}(zq^{m/2+1/2})
(x^+(z))(x^+(zq^{2}))\cdot\cdot\cdot
(x^+(zq^{(2m)})) \psi^{}(zq^{(m+1)3/2}), $$
$$ (x^-(zq^{2m+2}))(x^-(zq^{2m}))
\cdot\cdot\cdot
(x^-(zq^2))=$$
$$\vphi(zq^{-(m+1)3/2})(x^-(zq^{2m}))(x^-(zq^{2m-2}))
\cdot\cdot\cdot (x^-(z))\psi^{-1}(zq^{-m/2+1/2}),$$
and,  $(x^+(z))\cdot\cdot\cdot
(x^+(zq^{(2m)}))$ and $(x^-(zq^{(2m}))(x^-(zq^{2m-2}))
\cdot\cdot\cdot (x^-(z))$ are vertex operators. 
In the last section, 
we apply this method  to derive the parafermions $\phi_i(z)$, $i=0,...,m-1$,
on the module $\otimes^{m+1} \c F$  for  $U_q(\hat {\frak sl}(2))$,
such that 
$$(x^+(z))(x^+(zq^{2}))\cdot\cdot\cdot
(x^+(zq^{2k}))= V^+_k(z)\phi^{m,k}(z),$$
$$(x^-(zq^{2k}))(x^+(zq^{2k-2}))\cdot\cdot\cdot
(x^+(z))= V^-_k(z)\phi_{m-k-1}(z),$$
where $\phi_i(z)$ commute with $\vphi(z)$ and $\psi(z)$, and 
$ V^+_k(z)$ and $V^-_k(z)$ are vertex operators. 

\section{Quantum difference equation and vertex operator}
We will first present the current realization of $U_q(\hgtsl_2)$ given by 
Drinfeld\cite{Dr2}.

\begin{defn}
The algebra $U_q(\hgtsl_2)$ is an associative algebra with unit 
1 and the generators: $\vphi(m)$,$\psi(-m)$, $x^{\pm}(l)$, for 
$i=1,...,n-1$, $l\in \nz $ and $m\in \nz_{\leq 0}$ and a central
 element $c$. Let $z$ be a formal variable and 
 $x^{\pm}(z)=\sum_{l\in \nz}x^{\pm}(l)z^{-l}$, 
$\vphi(z)=\sum_{m\in \nz_{\leq 0}}\vphi(m)z^{-m}$ and 
$\psi(z)=\sum_{m\in \nz_{\geq 0}}\psi(m)z^{-m}$. In terms of the 
formal variables, 
the defining relations are 
\begin{align*}
&\vphi(0)\psi(0)=\psi(0)\vphi(0)=1, \\
& \vphi(z)\vphi(w)=\vphi(w)\vphi(z), \\
& \psi(z)\psi(w)=\psi(w)\psi(z), \\
& \vphi(z)\psi(w)\vphi(z)^{-1}\psi(w)^{-1}=
  \frac{g(z/wq^{-c})}{g(z/wq^{c})}, \\
& \vphi(z)x^{\pm}(w)\vphi(z)^{-1}=
  g(z/wq^{\mp \frac{1}{2}c})^{\pm1}x^{\pm}(w), \\
& \psi(z)x^{\pm}(w)\psi(z)^{-1}=
  g(w/zq^{\mp \frac{1}{2}c})^{\mp1}x^{\pm}(w), \\
& [x^+(z),x^-(w)]=\frac{1}{q-q^{-1}}
  \left\{ \delta(z/wq^{-c})\psi(wq^{\frac{1}{2}c})-
          \delta(z/wq^{c})\vphi(zq^{\frac{1}{2}c}) \right\}, \\
& (z-q^{\pm a}w)x^{\pm}(z)x^{\pm}(w)=
  (q^{\pm a}z-w)x^{\pm}(w)x^{\pm}(z), 
\end{align*}
where
\[ \delta(z)=\sum_{k\in \nz}z^k, \quad
   g(z)=\frac{q^{a}z-1}{z-q^{a}}\quad \text{(expanded around
 $z=0$) },
\quad a=2. \]
\end{defn}

For this current realization, 
Drinfeld also gave the Hopf algebra structure.

\begin{thm}
The algebra $U_q(\hgtsl_2)$ has a Hopf algebra structure, which are given 
by the following formulae. 

\noindent{\bf Coproduct $\Delta$}
\begin{align*}
\text{(0)}& \quad \Delta(q^c)=q^c\otimes q^c, \\
\text{(1)}& \quad \Delta(x^+(z))=x^+(z)\otimes 1+
            \vphi(zq^{\frac{c_1}{2}})\otimes x^+(zq^{c_1}), \\
\text{(2)}& \quad \Delta(x^-(z))=1\otimes x^-(z)+
            x^-(zq^{c_2})\otimes \psi(zq^{\frac{c_2}{2}}), \\
\text{(3)}& \quad \Delta(\vphi(z))=
            \vphi(zq^{-\frac{c_2}{2}})\otimes\vphi(zq^{\frac{c_1}{2}}), \\
\text{(4)}& \quad \Delta(\psi(z))=
            \psi(zq^{\frac{c_2}{2}})\otimes\psi(zq^{-\frac{c_1}{2}}),
\end{align*}
where $c_1$ is the action of the central element $c$ on the 
first component and $c_2$ is the action of the 
central element $c$ on the second component.  

\noindent{\bf Counit $\vep$}
\begin{align*}
\vep(q^c)=1 & \quad \vep(\vphi(z))=\vep(\psi(z))=1, \\
            & \quad \vep(x^{\pm}(z))=0.
\end{align*}
\noindent{\bf Antipode $\quad a$}
\begin{align*}
\text{(0)}& \quad a(q^c)=q^{-c}, \\
\text{(1)}& \quad a(x^+(z))=-\vphi(zq^{-\frac{c}{2}})^{-1}
                               x^+(zq^{-c}), \\
\text{(2)}& \quad a(x^-(z))=-x^-(zq^{-c})
                               \psi(zq^{-\frac{c}{2}})^{-1}, \\
\text{(3)}& \quad a(\vphi(z))=\vphi(z)^{-1}, \\
\text{(4)}& \quad a(\psi(z))=\psi(z)^{-1}.
\end{align*}

\end{thm}
This comultiplication structure requires certain completion
on the tensor space. For certain  representations, such 
as the $n$-dimensional representations of $U_q(\hgtsl_n)$ at a special value, 
this comultiplication may  not be  well-defined. 
Nevertheless, for 
any two highest weight representations, this comultiplication
is well-defined, because the action of the operator as a  coefficient  of 
$z^m$ of the currents operators on any element of such 
a module are zero if $m$ is small enough.

We will present the Frenkel-Jing construction of level $1$ representation 
of  $U_q(\hgtsl_2)$ on the Fock space. 

Consider an algebra generated by 
$\{ a_{k}|~k\in \nz \setminus \{ 0\} \}$ 
satisfying:
\[ [a_{k},a_{l}]=\delta_{k+l,0}\frac{[2k][k]}{k}, \]
where $[k]=\frac {q^k-q^{-k}}{q-q^{-1}}$. 
We call it the Heisenberg algebra. 

Let $ \overline{Q}=\nz \alpha $ be the root lattice of ${\frak sl}(2)$. 
Let us define a group algebra $\nc(q) [\overline{\cP}]$, where
 $\overline{\cP}$ is the weight lattice of $\gtsl_2$. 
Let ${\Lambda_1}$ be the fundamental weight of ${\frak sl}(2)$
 and $2{\Lambda_1}=a$. Let $\Lambda_0=0$. 

Set
\[ \cF_{i}:=\nc(q) [a_{-k}(~k\in \nz_{>0})]\otimes
            \nc(q) [\overline{Q}]
              e^{\overline{\Lambda}_i}
  . \]
This gives the Fock space. 

The action of operators 
$a_{k},\partial_{\alpha},e^{\alpha}$~$(1\leq j \leq N)$ is given by 
\begin{align*}
a_{k}\cdot f\otimes e^{\beta}& =\begin{cases}
                      a_{k}f \otimes e^{\beta}          & k< 0;  \\
            \text{$[a_{k},f]$}  \otimes e^{\beta} \quad & k> 0,
                                  \end{cases} \\
\partial_{\alpha}\cdot f\otimes e^{\beta}& 
        =(\alpha,\beta)f\otimes e^{\beta}
         \qquad \text{for}~f\otimes e^{\beta}\in \cF_{i}, \\
e^{\alpha}\cdot f\otimes e^{\beta} & 
=f\otimes e^{\alpha}e^{\beta}.
\end{align*}

\begin{lemma} The following 
action  on $\cF_{i}$ of 
$U_q(\hat {\frak sl}(2))$ gives a  level 1 highest weight  representation
 with  the $i$-th  fundamental
weight. 
\begin{align*}
\circ & \quad x^{\pm}(z)\mapsto
        \exp[\pm \sum_{k>0}\frac{a_{-k}}{[k]}q^{\mp \frac{1}{2}k}z^k]
        \exp[\mp \sum_{k>0}\frac{a_{k}}{[k]}q^{\mp \frac{1}{2}k}z^{-k}]
        e^{\pm \alpha}z^{\pm \partial_{\alpha}+1}, \\
\circ & \quad \vphi(z)\mapsto
        \exp[-(q-q^{-1})\sum_{k>0}a_{-k}z^k]q^{-\partial_{\alpha}}, \\
\circ & \quad \psi(z)\mapsto
        \exp[(q-q^{-1})\sum_{k>0}a_{k}z^{-k}]q^{\partial_{\alpha}}.
\end{align*}
\end{lemma}

This implies that on $\cF_{i}$ for the case of $U_q(\hgtsl_2)$
$$x^+(z)x^+(w)= z^2(1-\frac w z)(1-\frac {w}{zq^2}):x^+(z)x^+(w):$$
$$x^-(z)x^-(w)= z^2(1-\frac w z)(1-\frac {wq^2}{z}):x^-(z)x^-(w):$$
$$x^+(z)\vphi(w)= 
q^{-2}\frac {(\frac {wq^{-1/2}}z-q^2)}{(\frac {wq^{5/2}}z-1)}:\vphi(w)x^+(z):
=\frac {(\frac {wq^{-1/2}}z-q^2)}{(\frac {wq^{5/2}}z-1)}\vphi(w)x^+(z) $$
$$\psi(w)x^-(z)= 
\frac{(\frac {zq^{1/2}}w-q^2)} {(\frac {zq^{1/2}}w-q^2)}:\psi(w)x^-(z):
=\frac {(\frac {zq^{1/2}}w-q^2)}{(\frac {zq^{5/2}}w-1)}x^-(z)\psi(w)$$

\begin{lemma} 
Set $\c F=\Sigma\oplus \cF_{i}$. Any level $m$ integrable module 
is a submodule of  $\otimes^m \c F$. 
\end{lemma}

For the case of $\hat {\frak sl}(2)$, we have that the 
correlation functions of $e(z)e(w)$ and $f(z)f(w)$ have no poles, 
which are always polynomials of $
z,z^{-1},w,w^{-1}$. By the  correlation functions of an operator, we 
mean all the matrix coefficients of the operator. 
However, for the quantum case, we have \cite{DM}

\begin{thm}
For any level $m>1$ integrable module of $U_q(\hgtsl_2)$, 
the correlation functions of $x^+(z)x^+(w)$ has at most  poles  at 
$zq^{-2}=w$.  For any level $m>1$ integrable module of $U_q(\hgtsl_2)$, 
the correlation functions of $x^-(z)x^-(w)$ has at most poles  at 
$zq^{2}=w$. 
\end{thm}  

\begin{thm}
For any level $m$ integrable module of $U_q(\hgtsl_2)$, 
the correlation functions of $x^+(z_{m+1})x^+(z_{m})...x^+(z_2)x^+(z_1)$ 
is zero at $z_{i}/z_{i+1}=q^2$.For any level $m$ integrable module of $U_q(\hgtsl_2)$, 
the correlation functions of $x^-(z_{m+1})x^-(z_{m})...x^-(z_2)x^-(z_1)$ 
is zero, if $z_{i+1}/z_{i}=q^2$. 
\end{thm}

\begin{lemma}
On $\c F$,we have 
$$x^{\pm}(q^2z)=
q^{\pm \partial_{\alpha}}\exp[\pm \sum_{k>0}a_{-k}(q-q^{-1})q^k
q^{\mp \frac{1}{2}k}z^k]x^{\pm}(z)\times $$
$$q^{\pm \partial_{\alpha}}
\exp[\mp \sum_{k>0}(q-q^{-1})(-q^{-k})
{a_{k}}q^{\mp \frac{1}{2}k}z^{-k}]=$$
$$\vphi^{\mp 1}(zq^{\mp\frac{1}{2}+1})x^{\pm}(z)\psi^{\pm 1}(zq^{\pm
\frac{1}{2}+1}).$$
\end{lemma}

Proof. 
$$x^{\pm}(q^2z)=
        \exp[\pm \sum_{k>0}\frac{a_{-k}}{[k]}q^{\mp \frac{1}{2}k}q^{2k}z^k]
        \exp[\mp \sum_{k>0}\frac{a_{k}}{[k]}q^{\mp \frac{1}{2}k}q^{-2k}z^{-k}]
        e^{\pm \alpha}z^{\pm \partial_{\alpha}+1}q^{2\pm \partial_{\alpha}}
q^2= $$
$$\exp[\pm \sum_{k>0}\frac{a_{-k}}{[k]}q^{\mp \frac{1}{2}k}(q^{2k}-1)z^k]
x^{\pm}(z)
\exp[\mp \sum_{k>0}\frac{a_{k}}{[k]}q^{\mp \frac{1}{2}k}(q^{-2k}-1)z^{-k}]
q^2 q^{2\pm \partial_{\alpha}}=$$
$$\exp[\pm \sum_{k>0}a_{-k}(q-q^{-1})q^nq^{\mp \frac{1}{2}k}z^k]
x^{\pm}(z)\exp[\mp \sum_{k>0}(q-q^{-1})(-q^{-n})
{a_{k}}q^{\mp \frac{1}{2}k}z^{-k}]q^2 q^{2\pm \partial_{\alpha}}=
$$
$$q^{\pm \partial_{\alpha}}\exp[\pm \sum_{k>0}a_{-k}(q-q^{-1})q^k
q^{\mp \frac{1}{2}k}z^k]x^{\pm}(z)q^{\pm \partial_{\alpha}}
\exp[\mp \sum_{k>0}(q-q^{-1})(-q^{-k})
{a_{k}}q^{\mp \frac{1}{2}k}z^{-k}].$$

\begin{lemma}
On $\otimes^{m+1} \c F$, let $\Delta^m(x^+(z))=
(1\otimes 1\otimes ...\otimes \Delta)....(1\otimes \Delta)(\Delta)
((x^+)(z)),$ then 
$$\Delta^m(x^+(z))\Delta^m(x^+(zq^{2}))\cdot\cdot\cdot
\Delta^m(x^+(zq^{(2m)}))=$$
$$\vphi(zq^{1/2})\vphi(zq^{1/2+2})\cdot\cdot\cdot \vphi(zq^{1/2+2m-2})
x^+(zq^{2m})\otimes \vphi(zq^{3/2})\vphi(zq^{3/2+2})\cdot$$
$$\cdot\cdot
\vphi(zq^{3/2+2m-4})x^+(zq^{2m-1})\otimes 
 ......\otimes 
\vphi(zq^{m-1/2})x^+(zq^{m+1})\otimes  x^+(zq^m).$$
\end{lemma}

Proof. 
From the comultiplication formula,  on  $\otimes^{m+1} \c F$ 
we have 
$$(1\otimes 1\otimes ...\otimes \Delta)....(1\otimes \Delta)(\Delta)
((x^+)(z))=$$ 
$$x^+(z)\otimes 1....\otimes 1+ ...
+\vphi(zq^{1/2})\otimes 
\vphi(zq^{3/2})\otimes...\otimes x^+(zq^l)\otimes$$ 
$$ 1...\otimes 1
+....
+\vphi(zq^{1/2})\otimes\vphi(zq^{3/2})\otimes...\otimes
\vphi(zq^{m-1/2})\otimes x^+(zq^m)=$$
$$ \sum_{i=1}^{m} X^{+m}_{i}(z).$$

Let $$f^{m}_{a_1,...a_{m+1}}(z_1,....,z_{m+1})=X^{+m}_{a_1}(z_1)\times 
\cdot\cdot\cdot X^{+m}_{a_{m+1}}(z_{m+1}),$$
First, we know that there are  no poles at $z_i/z_{i+1}=q^2$ (Theorem 2.1). 
Let $0<i<j\leq m+1$
$$ (\vphi(zq^{1/2})\otimes 
\vphi(zq^{3/2})\otimes...\otimes x^+(zq^i)\otimes 1..\otimes 1)
(\vphi(wq^{1/2})\otimes 
\vphi(wq^{3/2})\otimes...\otimes x^+(wq^j)\otimes 1..\otimes 1)=$$
$$\vphi(zq^{1/2})\vphi(wq^{1/2})\otimes \vphi(zq^{3/2})\vphi(wq^{3/2})\otimes 
..$$
$$\otimes : x^+(zq^{i-1})\vphi(wq^{i-1/2}):\otimes \vphi(wq^{i+1/2})\otimes....
x^+(wq^j)\otimes ...\otimes 1 
\frac {(\frac {w}z-q^2)}{(\frac w{z}q^2)-1}.$$
Thus know that if there is any $i$, such
 that $a_i<a_{i+1}$, then the correlation functions of 
$f^m_{a_1,...a_{m+2}}(z_1,....,z_{m+2})$ are zero at $z_i/z_{i+1}=q^2$. 
That means the elements that are possibly not zero have the 
property that  $a_i\leq a_j$, if $i>j$. Because  
$x^+(z)x^+(zq^2)$ is zero, we have that 
the elements that are possibly not zero have the 
property that  $a_i< a_j$, if $i>j$. Therefore we finish the proof.

Let $x^{+m}(z)= \Delta^m(x^+(z))\Delta^m(x^+(zq^{2}))\cdot\cdot\cdot
\Delta^m(x^+(zq^{(2m)})).$

\begin{thm}
$$x^{+m}(zq^2)=\vphi^{-1}(zq^{1/2})\otimes \vphi^{-1}(zq^{3/2})\otimes....
\vphi^{- 1}(zq^{\frac{1}{2}+m})x^{+m}(z)$$
$$\psi (zq^{
\frac{1}{2}+2m+1})\otimes .....\otimes \psi^{}(zq^{
\frac{1}{2}+1+m})=$$
$$(1\otimes 1\otimes ...\otimes \Delta)....(1\otimes \Delta)(\Delta)
\vphi^{-1}(zq^{m/2+1/2})x^{+m}(z)$$
$$
(1\otimes 1\otimes ...\otimes \Delta)....(1\otimes \Delta)(\Delta)
 \psi^{}(zq^{(m+1)3/2}). $$
\end{thm}

Proof. 
$$x^{+m}(zq^2)= 
\vphi(zq^2q^{1/2})\vphi(zq^2q^{1/2+2})\cdot\cdot\cdot \vphi(zq^2q^{1/2+2m-2})$$
$$
x^+(zq^2q^{2m})\otimes \vphi(zq^2q^{3/2})\vphi(zq^2q^{3/2+2})\cdot\cdot\cdot
\vphi(zq^2q^{3/2+2m-4}x^+(zq^2q^{2m-1})\otimes $$
$$ ......\otimes 
\vphi(zq^2q^{m-1/2})x^+(zq^2q^{m+1})\otimes  x^+(zq^2q^m)=$$
$$
\vphi(zq^2q^{1/2})\vphi(zq^2q^{1/2+2})\cdot\cdot\cdot \vphi(zq^2q^{1/2+2m-2})
\vphi^{-1}(zq^{-\frac{1}{2}+1+2m})x^{+}(zq^{2m})\psi^(zq^{
\frac{1}{2}+1+2m})\otimes $$
$$
 \vphi(zq^2q^{3/2})\vphi(zq^2q^{3/2+2})\cdot\cdot\cdot
\vphi(zq^2q^{3/2+2m-4} 
\vphi^{-1}(zq^{-\frac{1}{2}+2m}x^{+}(zq^{2m-1})\psi^(zq^{
\frac{1}{2}+2m})\otimes ....$$ $$\otimes 
\vphi(zq^2q^{m-1/2}) \vphi^{-1}(zq^{-\frac{1}{2}+1+m+1})
x^{+}(zq^{m+1})\psi^{}(zq^{
\frac{1}{2}+1+m+1})\otimes \vphi^{-1}(zq^{-\frac{1}{2}+1+m})
\otimes $$ $$x^{+}(zq^{m})\psi^{}(zq^{
\frac{1}{2}+m+1})= 
\vphi^{-1}(zq^{1/2})\otimes \vphi^{-1}(zq^{3/2})\otimes....
\vphi^{- 1}(zq^{\frac{1}{2}+m})x^{+m}(z)\times $$
$$\psi(zq^{
\frac{1}{2}+2m+2})\otimes .....\otimes \psi^{}(zq^{\frac{1}{2}+1+m})
.$$

On  the module $\otimes^{m+1}{\c F}$, from the 
comultiplication formula, we have 
$$\Delta^m(a_{\pm n})= 
(1\otimes 1\otimes ...\otimes \Delta)....(1\otimes \Delta)(\Delta)a_{\pm n}
= $$
$$q^{(m)n/2} a_{\pm n} \otimes 1..\otimes 1 + 1\otimes 
q^{(m-2)n/2} a_{\pm n} \otimes 1..\otimes 1
+...+$$
$$1\otimes ....\otimes 1 \otimes q^{-(m)n/2} a_{\pm n}. $$
Then
$$[\Delta^m(a_{ n}), \Delta^m(a_{-n})]=[2n][(m+1)n]/n.$$
From the theorem above,  we have that 
$$x^{+m}(z)= \exp[\sum_{k>0}\frac{\Delta^m(a_{-k})}{[k]}q^{m/2-1/2} z^k]
        \exp[- \sum_{k>0}\frac{\Delta^m(a_{k})}{[k]}q^{-3m/2-1/2}
z^{-k}]\times$$
$$ (e^\alpha \otimes e^\alpha.....\otimes e^{\alpha})
(z^{\partial_{\alpha}+1}\otimes .....\otimes  z^{\partial_{\alpha}+1})
(q^{m\partial_{\alpha}}\otimes....\otimes q^{m\partial_{\alpha}})
q^{(1+m)m/2}. $$

\begin{cor}
$x^{+m}(z)$ is an vertex operator. 
\end{cor}

Similarly we can derive corresponding results for $x^-(z)$. 

\begin{lemma}
On $\otimes^{m+1} \c F$, let $\Delta^m(x^-(z))=
(1\otimes 1\otimes ...\otimes \Delta)....(1\otimes \Delta)(\Delta)
((x^-)(z)),$ then 
$$\Delta^m(x^-(zq^{(2m}))\Delta^m(x^-(zq^{2m-1}))\cdot\cdot\cdot
\Delta^m(x^-(z))=$$
$$x^-(zq^m)\otimes x^-(zq^{m+1})\psi(zq^{3/2+2m-2})\otimes ... 
\otimes $$
$$x^-(zq^{2m-1})\psi(zq^{3/2+2m-4})\cdot\cdot\cdot\psi(zq^{3/2+2})
\psi(zq^{3/2})\otimes x^-(zq^{2m})\psi(zq^{1/2+2m-2})\cdot\cdot\cdot
\psi(zq^{1/2+2})\psi(zq^{1/2}).$$
\end{lemma}

Let $ x^{-m}(z)=\Delta^m(x^-(zq^{(2m}))\Delta^m(x^-(zq^{2m-2}))
\cdot\cdot\cdot
\Delta^m(x^-(z)$.
\begin{thm}
$$x^{-m}(zq^2)=\vphi^{}(zq^{\frac{1}{2}+1+m})\otimes
...\otimes \vphi^(zq^{
\frac{1}{2}+2m+1}) x^{-m}(z)$$
$$\psi(z)^{- 1}(zq^{\frac{1}{2}+m}\otimes...
\psi^{-1}(zq^{3/2})\otimes \psi^{-1}(zq^{1/2})=$$
$$(1\otimes 1\otimes ...\otimes \Delta)....(1\otimes \Delta)(\Delta)
\vphi(zq^{-(m+1)3/2})x^{-m}(z)$$
$$(1\otimes 1\otimes ...\otimes \Delta)....(1\otimes \Delta)(\Delta)
\psi^{-1}(zq^{-m/2+1/2})$$
\end{thm}

Then, from the theorem above,  we have that 
$$x^{-m}(z)=\exp[- \sum_{k>0}\frac{\Delta^m(a_{-k})}{[k]}q^{(1+3m/2)k}z^k]
        \exp[ \sum_{k>0}\frac{\Delta^m(a_{k})}{[k]}q^{(1-m)/2)k}z^{-k}]
        $$
$$ (e^{-\alpha} \otimes e^{-\alpha}.....\otimes e^{-\alpha})
(z^{-\partial_{alpha}+1}\otimes .....\otimes  z^{-\partial_{alpha}+1})
(q^{-m\partial_{alpha}}\otimes....\otimes q^{-m\partial_{alpha}})
q^{(1+m)m/2}$$

\begin{cor}
$x^{-m}(z)$ is an vertex operator. 
\end{cor}


\section{Quantum parafermions}

In this section, we will derive quantum  parafermion  and 
explain the parafermionic construction of 
integrable modules of $U_q(\hat {\frak sl}(2))$ following the line 
of the work of Lepowsky and Wilson \cite{LW1} \cite{LW2}. This type 
of construction for the classical case
was also given in \cite{FZ} from a different point of 
view.

On the module $\otimes ^{m+1} \c F$, we have that 
$$(1\otimes 1\otimes ...\otimes \Delta)....(1\otimes \Delta)(\Delta)
((x^+)(z))=$$ 
$$x^+(z)\otimes 1....\otimes 1+ ...
+\vphi(zq^{1/2})\otimes 
\vphi(zq^{3/2})\otimes...\otimes x^+(zq^l)\otimes 1...\otimes 1$$ 
$$+....
+\vphi(zq^{1/2})\otimes\vphi(zq^{3/2})\otimes...\otimes
\vphi(zq^{m-1/2})\otimes x^+(zq^m).$$

So the $l+1$-th term   $X^{+m}_{l}$ is
$$
\vphi(zq^{1/2})\otimes 
\vphi(zq^{3/2})\otimes...\otimes \vphi(zq^{l-1/2})\otimes x^+(zq^l)
\otimes 1...\otimes 1= $$
$$ \exp[-(q-q^{-1})\sum_{k>0}a_{-k}q^{k/2}z^k]q^{-\partial_{\alpha}}
\otimes 
 \exp[-(q-q^{-1})\sum_{k>0}a_{-k}q^{3k/2}z^k]q^{-\partial_{\alpha}}\otimes 
$$
$$....\otimes 
 \exp[-(q-q^{-1})\sum_{k>0}a_{-k}q^{k(l-1/2)}z^k]q^{-\partial_{\alpha}}
\otimes 
 \exp[ \sum_{k>0}\frac{a_{-k}}{[k]}q^{- \frac{1}{2}k}q^{lk}z^k]
$$
$$        \exp[ -\sum_{k>0}\frac{a_{k}}{[k]}q^{- \frac{1}{2}k}q^{-lk}
z^{-k}]         e^{ \alpha}z^{ \partial_{\alpha}+1}q^{l \partial}q^l=
$$
$$
\exp[-(q-q^{-1})\sum_{k>0}z^k(a_{-k}q^{k/2}\otimes 1...\otimes 1+
1\otimes a_{-k}q^{3k/2}\otimes ....\otimes 1 +....+$$
$$1\otimes ...\otimes 1 \otimes a_{-k}q^{l-1/2}\otimes 1...\otimes 1)
+ 
1\otimes ..\otimes 1\otimes \frac{a_{-k}}{[k]}q^{- \frac{1}{2}k}q^{lk}\otimes 
1...\otimes 1]$$
$$ \exp[1\otimes ..\otimes 1
 -\sum_{k>0}\frac{a_{k}}{[k]}q^{- \frac{1}{2}k}q^{-lk}\otimes 1...\otimes 1]
(1\otimes ....\otimes 1\otimes  e^{ \alpha}\otimes 1....\otimes 1)
$$
$$q^{-\partial_{\alpha}}\otimes ....\otimes q^{-\partial_{\alpha}}
\otimes
 z^{ \partial_{\alpha}+1}q^{l \partial_{\alpha}}\otimes 1...\otimes 1)q^l .$$
 
\begin{lemma}
$$[\Delta^m(a_{k}), -(q-q^{-1})(a_{-k}q^{k/2}\otimes 1...\otimes 1+
1\otimes a_{-k}q^{3k/2}\otimes ....\otimes 1 +....+
1\otimes $$
$$...\otimes 1 \otimes a_{-k}q^{l-1/2}\otimes 1...\otimes 1)
+1\otimes ..\otimes 1\otimes \frac{a_{-k}}{[k]}q^{- \frac{1}{2}k}q^{lk}\otimes 
1...\otimes 1]=$$
$$ [2k]q^{-k(m+1)/2},$$

$$
[1\otimes ..\otimes 1
 -\sum_{k>0}\frac{a_{k}}{[k]}q^{- \frac{1}{2}k}q^{-lk}\otimes 1...\otimes 1, 
\Delta^m(a_{-k})]= $$
$$-q^{((-m-1)/2+l)k}[2k]/n (q^{- \frac{1}{2}k}q^{-lk})=$$
$$-q^{-k(m+1)/2)}[2k]/n. $$
\end{lemma}

Proof. 
$$[\Delta^m(a_{k}), -(q-q^{-1})(a_{-k}q^{k/2}\otimes 1...\otimes 1+
1\otimes a_{-k}q^{3k/2}\otimes ....\otimes 1 +....+
1\otimes $$
$$...\otimes 1 \otimes a_{-k}q^{l-1/2}\otimes 1...\otimes 1)
+1\otimes ..\otimes 1\otimes \frac{a_{-k}}{[k]}q^{- \frac{1}{2}k}q^{lk}\otimes 
1...\otimes 1]=$$
$$[2k][k]/n (q^{k/2}q^{k(-m)/2}+ q^{3/2k}q^{k(-m+2)/2}+..+
q^{(l-1/2)k}q^{((-m)/2+l-1)k})(-q+q^{-1})
+ $$
$$q^{-k/2}q^{lk}q^{((-m)/2+l)k}/[k])= [2k][k]/n(
q^{k/2}q^{k(-m)/2}(1+...+ q^{2(l-1)k})+q^{-k/2}q^{lk}q^{((-m)/2+l)k}/[k]) 
= $$
$$(-q+q^{-1})[2k][k]/n(q^{k/2}q^{k(-m)/2}(1-q^{2lk}/(1-q^{2k})
+q^{-k/2}q^{lk}q^{((-m)/2+l)k}/[k])=$$
$$[2k]q^{-k/2}q^{k(-m)/2}/n=[2k]q^{-k(m+1)/2}$$

From the comultiplication formula, we have 

$$\Delta^m(x^-(z))=
(1\otimes 1\otimes ...\otimes \Delta)....(1\otimes \Delta)(\Delta)
((x^-)(z))= $$
$$1\otimes ....\otimes 1\otimes x^-(z) +
1\otimes ....\otimes 1\otimes x^-(zq^1)\otimes 
\psi^{}(zq^{l-\frac{1}{2}})\otimes....+$$
$$...+ x^-(zq^{m})\otimes \psi^{}(zq^{m-\frac{1}{2}})\otimes ...
\otimes \psi^{}(zq^{\frac{1}{2}})=$$
$$\sum_{i=1}^{m} X_i^{-m}(z)$$

$$X_i^{-m}(z)=1\otimes ....\otimes 1\otimes x^-(zq^{i-1})\otimes 
\psi(zq^{i-1-\frac{1}{2}})=$$
$$1\otimes..\otimes \exp[- \sum_{k>0}\frac{a_{-k}}{[k]}q^{{(i-1)}k+
\frac{1}{2}k}z^k]
        \exp[ \sum_{k>0}\frac{a_{k}}{[k]}q^{-{(i-1)}k+\frac{1}{2}k}z^{-k}]
        e^{- \alpha}z^{- \partial_{\alpha}+1}q^{-{(i-1)}\partial_{\alpha}}
q^{(i-1)} \otimes $$
$$q^{\partial_{\alpha}}
\exp[(q-q^{-1}\sum_{k>0}{a_{k}}q^{{(i-1)}k-\frac{1}{2}k}z^{-k}]\otimes ...
\otimes q^{\partial_{\alpha}}
\exp[(q-q^{-1}\sum_{k>0}{a_{k}}q^{\frac{1}{2}k}z^{-k}]=$$
 $$ \exp [1\otimes..\otimes- \sum_{k>0}\frac{a_{-k}}{[k]}q^{{(i-1)}k+
\frac{1}{2}k}z^k]\otimes 1..\otimes 1]\times $$
$$\exp[1\otimes 1...\otimes 1 
\sum_{k>0}\frac{a_{k}}{[k]}q^{-{(i-1)}k+\frac{1}{2}k}z^{-k}
\otimes 1...\otimes 1+...
1\otimes ..1\otimes 
(q-q^{-1})\sum_{k>0}{a_{k}}q^{{(i-1)}k-\frac{1}{2}k}z^{-k}\otimes
$$ $$1...\otimes 1
...+...+1\otimes ...\otimes 1
(q-q^{-1})\sum_{k>0}{a_{k}}q^{\frac{1}{2}k}z^{-k}\otimes 1..
\otimes 1]\times $$
$$(1\otimes ...\otimes 1\otimes  e^{- \alpha}\otimes 1...\otimes 1)(
1\otimes..1\otimes z^{- \partial_{\alpha}+1}q^{-{(i-1)}\partial_{\alpha}}
\otimes q^{\partial_{\alpha}} \otimes ...\otimes 
q^{\partial_{\alpha}}q^{(i-1)}.
$$

\begin{lemma}
$$
[\Delta^m(a_{k}),
1\otimes..\otimes-\frac{a_{-k}}{[k]}q^{{(i-1)}k+
\frac{1}{2}k}]\otimes 1..\otimes 1]= 
-[2k]/n(q^{mk/2+k/2})$$
$$[1\otimes 1...\otimes 1 
\frac{a_{k}}{[k]}q^{-{(i-1)}k+\frac{1}{2}k}z^{-k}
\otimes 1...\otimes 1+1\otimes ..1\otimes 
(q-q^{-1}){a_{k}}q^{{(i-1)}k-\frac{1}{2}k}z^{-k}\otimes
$$ 
$$1...\otimes 1
...+...+1\otimes ...\otimes 1
(q-q^{-1})\sum_{k>0}{a_{k}}q^{\frac{1}{2}k}z^{-k}\otimes 1..
\otimes 1, \Delta^m(a_{-k})]=$$
$$[2k]/n(q^{mk/2+k/2})$$
\end{lemma}

Let $$V(\pm m,z)= \exp [\frac {\Delta^m(a_{-k})}{[(m+1)k]}q^{\mp (m+1)k/2} z^k]
\exp [\frac {\Delta^m(a_{k})}{[(m+1)k]}q^{\mp (m+1)k/2} z^{-k}]\times$$
$$(e^{\pm\alpha/(m+1)}\otimes .....\otimes e^{\pm\alpha/(m+1)})
(z^{\pm\partial_\alpha/(m+1)}\otimes .....\otimes z^{\pm\partial\alpha/(m+1)})
z$$

\begin{lemma}
 $$V(\pm m,z)V(\pm m,w)=
 \exp [-\frac {2k}{k[(m+1)k]}q^{\mp(m+1)k} (z/w)^{-k}]z^{2/(m+1)}
:V(\pm m,z)V(\pm m,w):$$

$$V(\pm m,z)V(\mp m,w)= 
 \exp [\frac {2k}{k[(m+1)k]} 
(z/w)^{-k}]z^{-2/(m+1)}e^{-2\pi i/(m+1)}
:V(\pm m,z)V(\mp m,w):$$
\end{lemma}

We will denote 
$\exp [-\frac {2k}{[(m+1)k]}q^{\mp(m+1)k} (z/w)^{-k}]z^{2/(m+1)
}$  by $f^{\pm}(w,z)$  and 
$\exp [\frac {2k}{[(m+1)k]} (z/w)^{-k}]z^{-2/(m+1)}e^{-2 
\pi i/(m+1)}
$ by 
 $p(w,z)$. 

Let $X^{\pm m}(z)= V(\pm m,z)\phi^{\pm m}(z)$ and 
$X_i^{\pm m}(z)= V(\pm m,z)\phi_i^{\pm m}(z)$. We have that 
$\phi^{\pm m}(z)=\sum \phi_i^{\pm m}(z)$. 
\begin{prop}
On the space $\otimes^{m+1}\c F$,
$$\phi^{\pm m}(z)\Delta^m(\vphi(w))= 
\Delta^m(\vphi(w))\phi_i^{\pm m}(z),$$
$$\phi^{\pm m}(z)\Delta^m(\vphi(w))= 
\Delta^m(\vphi(w))\phi_i^{\pm m}(z),$$
$$\phi^{\pm m}(z)\Delta^m(\psi(w))= 
\Delta^m(\psi(w))\phi^{\pm m}(z),$$
$$\phi^{\pm m}(z)\Delta^m(\psi(w))= 
\Delta^m(\psi(w))\phi^{\pm m}(z).$$
\end{prop}

As the commutants  to $\Delta^m(\vphi(w)$ and $\Delta^m(\psi(w)$, 
$\phi_i^{\pm m}(z)$ degenerate into the classical parafermions 
respectively \cite{LW1}\cite{LW2} \cite{FZ}. Thus 
$\phi_i^{\pm m}(z)$ gives us the quantum parafermions. 

\begin{prop}
On the space $\otimes^{m+1}\c F$,
$$:\phi_{m+1}^{+m}(z)\phi_{m}^{+m}(zq^2)...
\phi_{1}^{+m}(zq^{2m}):=1;$$
$$:\phi_{1}^{-m}(z)\phi_{2}^{-m}(zq^{-2})...
\phi_{m+1}^{-m}(zq^{-2m}):=1.$$
\end{prop}
This is a quantum version of a classical relations\cite{LW2}\cite{FZ}, which 
comply with the results in \cite{DM}. 

From the calculation above, we can easily write down  the commutation 
relations between $\phi_i^{\pm m}(z)$ and $\phi_i^{\pm m}(w)$
 and the commutation relations between   
$\phi_i^{\pm m}(z)$ and $\phi_i^{\mp}(w)$. 

\begin{lemma}
If $i<j$,
$$X^{+m}_i(z)X^{+m}_j(w)= \frac {(\frac wz-q^2)}
{(\frac {wq^{2}}z-1)}X^{+m}_j(w)X^{+m}_i(z)=
g(w/z)^{-1}X^{+m}_j(w)X^{+m}_i(z),$$
$$X^{+m}_j(z)X^{+m}_i(w)= \frac {(\frac {zq^{2}}w-1)}{(\frac zw-q^2)}
X^{+m}_i(w)X^{+m}_j(z)=g(z/w)X^{+m}_i(w)X^{+m}_j(z),$$
$$X^{-m}_j(z)X^{-m}_i(w)=  \frac  {(\frac {zq^2}w-1)}
{\frac {z}w-q^2)}
X^{-m}_i(w)X^{-m}_j(z)=g(z/w)X^{-m}_i(w)X^{-m}_j(z) ,$$
$$X^{-m}_i(z)X^{-m}_j(w)= \frac {(\frac wz-q^2)}
 {(\frac {wq^{2}}z-1)}
X^{-m}_j(w)X^{-m}_i(z)=g(w/z)^{-1}X^{-m}_j(w)X^{-m}_i(z),$$
$$X^{\pm m}_i(z)X^{\mp m}_j(w)=X^{\mp m}_j(w)X^{\pm m}_i(z);$$
$$X^{+m}_i(z)X^{+m}_i(w)=z^2 \frac {(1-\frac wz)}
{(1-\frac {w}{zq^2})}:X^{+m}_i(w)X^{+m}_i(z):,$$
$$X^{-m}_i(z)X^{-m}_i(w)= z^2(1-\frac w z)(1-\frac {wq^2}{z})
:X^{-m}_i(z)X^{-m}_i(w):.$$
$$X^{+m}_i(z)X^{-m}_i(w)=z^{-2} \frac 1 {(1-\frac w{zq})
(1-\frac{w}{zq^{-1}})}:X^{-m}_i(w)X^{+m}_i(z):,$$
$$X^{-m}_i(z)X^{+m}_i(w)=w^{-2} \frac 1 {(1-\frac w{zq})
(1-\frac {w}{zq^{-1}})}:X^{-m}_i(z)X^{+m}_i(w):.$$
\end{lemma}

\begin{lemma}
If $i<j$,
$$f^+(w,z)\phi^{+m}_i(z)\phi^{+m}_j(w)=f^+(z,w) 
\frac {(\frac wz-q^2)}
{(\frac {wq^{2}}z-1)}\phi^{+m}_j(w)\phi^{+m}_i(z)=
$$ 
$$f^+(z,w)
g(w/z)^{-1}\phi^{+m}_j(w)\phi^{+m}_i(z),$$
$$f^+(w,z)\phi^{+m}_j(z)\phi^{+m}_i(w)=f^+(z,w) \frac 
{(\frac {zq^{2}}w-1)}{(\frac zw-q^2)}
\phi^{+m}_i(w)\phi^{+m}_j(z)=$$
$$f^+(z,w)g(z/w)\phi^{+m}_i(w)\phi^{+m}_j(z),$$
$$f^-(w,z)\phi^{-m}_j(z)\phi^{-m}_i(w)=f^-(z,w)
  \frac  {(\frac {zq^2}w-1)}
{\frac {z}w-q^2)}
\phi^{-m}_i(w)\phi^{-m}_j(z)=
$$
$$f^-(z,w)g(z/w)\phi^{-m}_i(w)\phi^{-m}_j(z) ,$$
$$f^-(w,z)\phi^{-m}_i(z)\phi^{-m}_j(w)=f^-(z,w) \frac {(\frac wz-q^2)}
 {(\frac {wq^{2}}z-1)}
\phi^{-m}_j(w)\phi^{-m}_i(z)=$$
$$f^-(z,w)g(w/z)^{-1}
\phi^{-m}_j(w)\phi^{-m}_i(z),$$
$$p^\pm(w,z)\phi^{\pm m}_i(z)\phi^{\mp m}_j(w)=
p^\pm(z,w)\phi^{\mp m}_j(w)\phi^{\pm m}_i(z);$$
$$\phi^{+m}_i(z)\phi^{+m}_i(w)=f^+(w,z)z^2 \frac {(1-\frac wz)}
{(\frac 1-{w}{zq^2}-1)}:\phi^{+m}_i(w)\phi^{+m}_i(z):,$$
$$f^+(w,z)\phi^{-m}_i(z)\phi^{-m}_i(w)=f^-(w,z) 
z^2(1-\frac w z)(1-\frac {wq^2}{z})
:\phi^{-m}_i(z)\phi^{-m}_i(w):.$$
$$\phi^{+m}_i(z)\phi^{-m}_i(w)=p(w,z)z^{-2} \frac 1 {(1-\frac w{zq})
(\frac 1-{w}{zq^{-1}})}:\phi^{-m}_i(w)\phi^{+m}_i(z):,$$
$$\phi^{-m}_i(z)\phi^{+m}_i(w)=p(w,z) w^{-2} \frac 1 {(1-\frac w{zq})
(\frac 1-{w}{zq^{-1}})}:\phi^{-m}_i(z)\phi^{+m}_i(w):.$$
\end{lemma}

\begin{thm}
$$
[p(w,z)\phi^{+m}(z)\phi^{-m}(w)-
p(z,w)\phi^{-m}(w)\phi^{+m}(z)]= \frac 1{q-q^{-1}}
(\delta(z/wq^{m+1})-(\delta(z/wq^{-(m+1)}),$$
$$f^\pm(w,z)(z-wq^{\pm2})\phi^{\pm(z)m}\phi^{\pm m}(w)
=f^\mp(z,w)(zq^{\pm2}-w)\phi^{\pm m}(w)\phi^{\pm m}(z).$$
\end{thm}

We will call $\phi^{\pm m}(z)$ quantum parafermions, which degenerates into 
the classical parafermions when $q$ goes to 1. 

\begin{cor}
$$\phi^{+1}(z)=\phi^{-1}(z),$$
$$\phi^{+1}_i(z)=\phi^{-1}_i(z), i=1,2;$$
$$\phi^{\pm 1}_i(z)\phi^{\pm}_i(w)=-
\phi^{\pm 1}_i(w)\phi^{\pm 1}_i(z), i=1,2;$$
$$\{\phi^{+1}_i(z), \phi^{-1}_i(w)\}= 
\delta(q^{-2}z/w). $$
\end{cor}
We can see that the quantum fermions is basically the same
as the classical ones, but with certain shifts.

\begin{lemma}
On the module $\otimes ^{m+1}\c F$, $N<m$, 
$$\Delta^m(x^+(z))\Delta^m(x^+(zq^2))
...\Delta^m(x^+(zq^{2N-2}))\Delta^m(x^+(zq^{2N}))= $$
$$\sum_{i_1>i_2....>i_{N+1}}X^{+m}_{i_1}(z)...X^{+m}_{i_{N+1}}(zq^{2N});$$
$$\Delta^m(x^-(zq^{2N}))\Delta^m(x^-(zq^{2N-2}))...
\Delta^m(x^-(zq^{2}))\Delta^m(x^-(z))= $$
$$\sum_{j_1<j_2....<j_{N+1}}X^{-m}_{j_1}(zq^{2n})...X^{-m}_{j_{N+1}}(z).$$
\end{lemma}
We can prove it with the same method as in Lemma 2.7.

\begin{lemma}
Let $ 0\leq N< m$,
$i_1>i_2....>i_{N+1}$, $j_1<j_2....<j_{m-N}$ and the set 
$\{ i_1,....,i_{N+1}, j_1,....,j_{m-N} \}$ is the set $
\{1,2,3,...,m+1\}$. On the module $\otimes ^{m+1}\c F$, 
let 
$$X^{+m}_{i_1}(z)...X^{+m}_{i_{N+1}}(zq^{2N})= $$
$$:V(+m,z)V(+m,zq^2)...V(+m,
zq^{2N}):\phi^{+m}_{i_1>i_2....>i_{N+1}}(z),$$
$$X^{-m}_{j_1}(zq^{2(m-N-1)})...X^{-m}_{j_{m-N+1}}(z)=
$$ $$:V(-m,z)V(-m,zq^2)...V(-m,
zq^{2(m-N-1)}):\phi^{-m}_{j_1<j_2....<j_{m-N+1}}(z).$$
Then $$\phi^{+m}_{i_1>i_2....>i_{N+1}}=\phi^{-m}_{j_1<j_2....<j_{m-N+1}}(z)$$
\end{lemma}

Proof. 
Let $$\bar x(z)= \exp[\sum_{k>0}\frac{a_{-k}}{[k]}q^{\frac{1}{2}k}z^k]
        \exp[- \sum_{k>0}\frac{a_{k}}{[k]}q^{\frac{1}{2}k}z^{-k}]
        e^{ \alpha}z^{\partial_{\alpha}-1}. $$
Then 
$$\bar x(zq^2)= \vphi(zq^{3/2})^{-1} \bar x(z)\psi(zq^{1/2}),$$
$$\bar x(zq^{-1})= x^+(z) \psi(zq^{-1/2})^{-1}q.$$

$$  \vphi(zq^{1/2+l})\vphi(zq^{1/2+2+l})\cdot\cdot\cdot \vphi(zq^{1/2+2m-l-2})
x^+(zq^{2m-l})=$$ $$ \vphi(zq^{1/2+l})\vphi(zq^{1/2+2+l})..
\vphi(zq^{1/2+2m-l-2})
x^+(zq^{2m})\psi(zq^{(2m-l)+3/2})^{-1}....\psi(zq^{2m+3/2})^{-1}.$$
$$ \vphi(zq^{1/2+l})\vphi(zq^{1/2+2+l})\cdot\cdot\cdot \vphi(zq^{1/2+2m-l-2})
x^+(zq^{2m-l})=$$ $$
\vphi(zq^{1/2+l})\vphi(zq^{1/2+2+l})\cdot\cdot\cdot \vphi(zq^{1/2+2m-l-2})
\bar x(zq^{2m-l-1})\psi(zq^{2m-l-1/2})q=$$ 
$$\vphi(zq^{1/2+l})\vphi(zq^{1/2+2+l})\cdot\cdot\cdot \vphi(zq^{1/2+2m-l-2k-2})
\bar x(zq^{2m-l-2k-1})\psi(zq^{2m-l-2k-1/2})\cdot$$ $$
\cdot\cdot \psi(zq^{2m-l-1-3/2})\psi(zq^{2m-l-1/2})q. $$
This gives us the proof.

\begin{cor}
For $0\leq N<m$,
let $$\phi^{+m,N}(z)=\sum \phi^{+m}_{i_1>i_2....>i_{N+1}}(z)$$ 
$$\phi^{-m,m-N-1}(z)=\sum \phi^{-m}_{i_1<i_2....<i_{m-N+1}}(z).$$
Then $$\phi^{+m,N}(z)=\phi^{-m,m-N-1}(z)$$. 
\end{cor}

Let $ V^-(\pm m,z)$ be an vertex operators such that 
$: V^-(\pm m,z) V(\pm m,z):=1$. 

\begin{cor}
$$: V^-(+ m,z)...V^-(+m,zq^{2N}):\Delta^m(x^+(z))....\Delta^m(x^+(zq{^2N}))=
$$
$$:V^-(- m,zq^{2m-2N})...V^-(-m,z):\Delta^m(x^-(zq^{2m-2N}))
....\Delta^m(x^-(z)).$$
\end{cor}
This generalizes the corresponding results in \cite{LP}.

\begin{defn}
Let $$\phi^{m,n}=\sum \phi^{+m}_{i_1>j_2...>i_{N+1}}(z),$$
for $N<m$. We call $\phi^{m,N}(z)$ quantum parafermions. 
\end{defn}
These quantum parafermions degenerate into the 
classical parafermions, when q goes to 1. 

With the results above, we can write the operator product expansion 
of these quantum parafermions and
the commutation relations between
those quantum parafermions. The operator product expansion is 
very complicated as it can been see from the commutation relations
below.

\begin{prop}
$$\prod_{j=0,1..N}^{i=0,1..,M}
(zq^{2j}-q^2wq^{2i})f^+(wq^{2i},zq^{2j})\phi^{m,N}(z)\phi^{m,M}(w)
=$$
$$\prod_{j=0,1..N}^{i=0,1..,M}
(q^2zq^{2j}-wq^{2i})f^+(zq^{2j},wq^{2i})\phi^{m,M}(w)\phi^{m,N}(z);$$
$$\prod_{j=0,1..m-N-1}^{i=0,1..,m-M-1}
(zq^{2j}-q^{-2}wq^{2i})f^-(wq^{2i},zq^{2j})\phi^{m,N}(z)\phi^{m,M}(w)
=$$
$$\prod_{j=0,1..m-N+1}^{i=0,1..,m-M+1}
(q^{-2}zq^{2j}-wq^{2i})f^-(zq^{2j},wq^{2i})\phi^{m,M}(w)\phi^{m,N}(z);$$
\end{prop}

A version of quantum parafermion is given in \cite{BV}, where
only $\phi^{m,0}(z)$ and $\phi^{m,m-1}(z)$ are given. 
With the Drinfeld comultiplication, we are able to follow 
the line of \cite{LW1} \cite{LW2} to derive all the 
integrable representation. This automatically leads us to the quantum 
parafermions characterized  as the commutant to $\vphi(z)$ and $\psi(z)$. 
Clearly, we can use the operators $\phi^{\pm m}(z)$ to derive $x^\pm(z)$, 
by $V(\pm,z)\phi^\pm(z)$, which gives parafermion constructions.
Combining the results in \cite{DM}, we bascially derive all the 
corresponding structure corresponding to the results in \cite{LW1} \cite{LW2}
concerning the structure of standard modules for the case of 
$\hat {\frak sl}(2)$, which essentially  prepares all the 
necessary conditions for the extension of the 
results in \cite{FS} to the quantum cases. 
On the other hand, parafermions in the classical case 
have very important applications in the  conformal field theory
\cite{FZ}. Our $\phi^{ m,N}(z)$ as deformed fermions should play an
important role in formulating the theory of the 
 quantization of the  conformal field  theory. The results in this
paper can be extended to the other cases of 
quantum affine algebras in a
straightforward way, which will gives us the structure of the standard 
modules of the quantum  affine algebras like in  \cite{LW1} \cite{LW2}.  

\bigskip
{\it Acknowledgement}

\medskip
\noindent
The authors would like to thank Tetsuji  Miwa for useful discussions.
J.D. is supported by the grant Reward research (A) 08740020 from the 
Ministry of Education of Japan. B.F. would like to thank 
RIMS for the hospitality throughout  his stay, during which 
this work was  done.

\end{document}